\documentclass[11pt,a4paper,reqno]{amsart}
\usepackage{amsmath,amssymb,amsfonts,amsthm,mathtools,mathrsfs}
\usepackage[T1]{fontenc}
\usepackage{lmodern,microtype}
\usepackage{enumerate,booktabs,array}
\usepackage[colorlinks=true,linkcolor=black,citecolor=black,urlcolor=blue]{hyperref}
\numberwithin{equation}{section}
\allowdisplaybreaks[2]
\newcommand{\R}{\mathbb R}
\newcommand{\Z}{\mathbb Z}

\newcommand{\dd}{\,\mathrm d}
\newcommand{\supp}{\operatorname{supp}}

\newcommand{\ip}[2]{\langle #1,#2\rangle}

\newcommand{\Ree}{\operatorname{Re}}

\makeatletter
\@namedef{subjclassname@2020}{\textup{2020} Mathematics Subject Classification}
\makeatother
\newtheorem{thm}{Theorem}[section]

\newtheorem{proposition}[thm]{Proposition}

\newtheorem{lemma}[thm]{Lemma}
\newtheorem{coro}[thm]{Corollary}
\newtheorem{corollary}[thm]{Corollary}

\newtheorem{ques}[thm]{Question}
\newtheorem{exam}[thm]{Example}
\theoremstyle{remark}

\newtheorem{remark}[thm]{Remark}
\title[On the spectrality of additive measures]
{Spectrality of additive measures with Bernoulli convolutions}
\author{Yi-Qiu Mao}
\address{Y.-Q. Mao, School of Mathematics and Information Science, Guangzhou University,
 Guangzhou, 510006, P.~R.~China}
\email{yqmao@gzhu.edu.cn}
\author{Zhi-Yi Wu}
\address{Z.-Y. Wu, School of Mathematics and Information Science, Guangzhou University,
 Guangzhou, 510006, P.~R.~China}
\email{zhiyi\_wu2021@163.com}
\begin{document}
\begin{abstract}
The spectral properties of additive measures of Lebesgue type are by now
completely settled \cite{LLP21,ALZ23,KW25,Lu25}. In this paper we continue
the study of spectrality for Cantor-type additive measures. For $0<r<1$, let
$\mu_r$ denote the Bernoulli convolution generated by the maps
$x\mapsto rx$ and $x\mapsto r(x+2)$ with equal weights, and consider the
Bernoulli-type additive measure
\[
 \mu_{r,s}=\tfrac12(\mu_r\times\delta_0+\delta_0\times\mu_s).
\]
For every even integer $R=2k$ we construct an explicit spectrum of the one-dimensional measure $\mu_{1/R}$, and from it an explicit spectrum of $\mu^{1/R}:=\mu_{1/R,1/R}$ lying on the antidiagonal is constructed. Combined with the recently established absence of Fourier frames for the Bernoulli
convolutions $\mu_{1/N}$ for odd integers $N\geq3$ \cite{DLT26,FSW26}, this yields a parity dichotomy: $\mu^{1/R}$ is spectral if and only if $R$ is even. For unequal contraction ratios, we prove that $\mu_{1/2,1/4}$ is not spectral. Several open problems and future directions are discussed.
\end{abstract}
\maketitle

\section{Introduction}\label{sec:intro}

A compactly supported Borel probability measure $\sigma$ on $\R^d$ is called a
{\it spectral measure} if it admits a countable set $\Lambda\subseteq\R^d$ for which
$E(\Lambda)=\{e_\lambda(x)=e^{2\pi i\langle\lambda,x\rangle}:\lambda\in\Lambda\}$ is
an orthonormal basis for $L^2(\sigma)$. Such a set $\Lambda$ is called  a {\it spectrum} of $\sigma$.

The study of spectral measures has its origins in the famous spectral set conjecture raised by Fuglede in 1974, which asserted that the normalised Lebesgue measure restricted on a bounded measurable set \(\Omega\) is spectral if and only if $\Omega$ tiles \(\mathbb{R}^d\) by translations \cite{Fug74}. Although the conjecture is now known to fail in general \cite{Tao04}, it has led to the development of many new research problems. Among these, in 1998, Jorgensen and Pedersen~\cite{JP98} discovered the first example of a non-atomic singular spectral measure, i.e., the Bernoulli convolution $\mu_r$ for $r=1/(2k)$, where $k\ge 2$ is an integer, and one of its spectra is as follows:
 \[
\Lambda(2k, C)=\left\{ \sum_{j=0}^{m} (2k)^j c_j : c_j \in \{0,k/2\}, \ m \ge 0 \right\}
\]
Since then, the spectrality of the Bernoulli convolutions has been studied extensively. Building on earlier work of Dutkay and Jorgensen \cite{DJ07}, Hu and Lau \cite{HL08}, Dai eventually settled the spectrality problem for Bernoulli convolutions, proving that $\mu_r$ is spectral if and only if $r=1/(2k)$ for some integer $k\ge1$ \cite{Dai12}. Meanwhile, the spectrality of self-similar measures \cite{AL23,AW21,DHLau14,DL15,LW02}, self-affine measures \cite{DFY24,DHLai19} and Moran measures has been extensively studied \cite{AFL19,AH14,AHLau15,DHLY24,LMW24,LLZ26,Shi19}. Also, various new phenomena different from those in the spectral theory for the Lebesgue measure were discovered. These mainly include:
\begin{itemize}
    \item the spectral eigenvalue problem (the spectral eigen-matrix
    problem in higher dimensions): for a given spectrum $\Lambda$ of a singular continuous spectral
    measure, characterizing those scalars $\lambda$ (or matrices $A$ in higher
    dimensions) such that $\lambda\Lambda$ (respectively $A\Lambda$) is again a
    spectrum \cite{ADH22,Dai16,DH16,HTW19}. For several classes of self-similar spectral measures in one dimension, this spectral eigenvalue problem has been studied \cite{Str00,Dai16,HTW19}. For example, Strichartz first discovered that $5\Lambda(4,\{0,1\})$ is a spectrum of $\mu_{1/4}$ \cite{Str00}. The problem is intimately linked to number-theoretic questions \cite{DH16}, and its higher-dimensional eigen-matrix counterpart has been investigated for the planar Sierpinski measures \cite{ADH22}.
    \item the distribution of the Beurling dimension of spectra. Dutkay et al. \cite{DHSW11} and He et al. \cite{HKTW18} proved that the Beurling dimension of the spectra of general self-similar spectral measure which satisfies the open set condition is not greater than the Hausdorff dimension of its support, and they may equal under some mild condition. Spectra of zero Beurling dimension and spectra of arbitrary intermediate Beurling dimension for some fractal measures were later constructed in \cite{AL23-0,DHLai13,LW24,Lu26}.
    \item the convergence of the associated mock Fourier series. Strichartz \cite{Str06} initiated the convergence theory of mock Fourier
series: for the Bernoulli convolutions \(\mu_{1/(2k)}\) with \(k\ge2\) and
the spectra \(\Lambda(2k,\{0,k/2\})\), he proved that the mock Fourier series
of every \(f\in C(\supp(\mu_{1/(2k)}))\) converges uniformly on
\(\supp(\mu_{1/(2k)})\); in the same paper he also studied the case \(k\ge3\)
and established pointwise convergence results in that setting. On the other
hand, Dutkay, Han and Sun \cite{DHS14} showed that convergence can fail for
other spectra: there exists \(f\in C([0,1])\) whose mock Fourier series with
respect to the spectrum \(17\Lambda(4,\{0,1\})\) diverges at \(0\).
Recently, Tang and the second author \cite{TW26} proved that, for the
spectrum \(\Lambda(4,\{0,1\})\), the mock Fourier series of every
\(f\in L^1(\mu_{1/4})\) converges to \(f\) \(\mu_{1/4}\)-almost everywhere
on \(\supp(\mu_{1/4})\); this complements Strichartz's pointwise
convergence results for \(k\ge3\) with the case \(k=2\), and solves Strichartz's long-standing problem \cite[p.~341]{Str06}.
\end{itemize}
For an excellent survey on this subject in one dimension, we refer the reader to
\cite{DFH26}.

Lai, Liu and Prince~\cite{LLP21} undertook a systematic study of additive measures obtained by placing two measures on the coordinate axes. This advanced the line of research initiated by Lev~\cite{Lev18} from the existence of Fourier frames for singular measures (i.e., the frame spectrality) to the exact characterization of spectrality for this specific additive model. For two non-atomic Borel probability measures $\mu$ and $\nu$, their
associated additive measure is
\[
 \frac12(\mu\times\delta_0+\delta_0\times\nu).
\]
Since \(\mu\) and \(\nu\) are non-atomic, this additive measure is continuous. Moreover, its support is concentrated on the union of the two coordinate axes, which has planar Lebesgue measure zero; hence, it is singular with respect to the planar Lebesgue measure. For the Lebesgue-type family
\[
 \frac12\big(\mathcal L_{[t_1,t_1+1]}\times\delta_0
       +\delta_0\times\mathcal L_{[t_2,t_2+1]}\big),
\]
where \(\mathcal L_{[a,b]}\) denotes the Lebesgue measure on
the interval \([a,b]\), the spectrality problem has been completely settled  \cite{LLP21,ALZ23,KW25,Lu25}: it is spectral if and only if $t_1+t_2\in\Z\setminus\{-1\}$ or $t_1-t_2\in\Z\setminus\{0\}$. In this paper, we consider the case where both component measures are the Bernoulli convolutions, with possibly distinct contraction ratios.

 We recall the definition of the Bernoulli convolution. For $r\in(0,1)$, the {\it Bernoulli convolution} $\mu_r$ is the probability measure on $\R$ that is the  distribution of the random variable \(\sum_{j=1}^\infty \xi_j r^j\), where \(\xi_i \) are i.i.d. random variables assuming values $0$ and $2$ with equal probability. Equivalently, $\mu_r$ is the unique Borel probability measure satisfying
\begin{equation}\label{eq:ifs}
\mu_r(\cdot)=\frac12\,\mu_r(\tau_1^{-1}(\cdot))+\frac12\,\mu_r(\tau_2^{-1}(\cdot)),
\end{equation}
where $\tau_1(x)=rx$ and $\tau_2(x)=r(x+2)$ (see \cite{Fal14}). This measure is supported on a compact $T_r$ as follows:
 \begin{equation}\label{eq:expand}
 T_r=\left\{\sum_{j=1}^\infty c_jr^j:~c_j\in\{0,2\}~\text{for all}~j\right\}.
 \end{equation}

We refer to the associated additive measure
\begin{equation}\label{eq:two-param}
\mu_{r,s}=\frac12(\mu_r\times\delta_0+\delta_0\times\mu_s)
\end{equation}
as the Bernoulli-type additive measure.

Although the one-dimensional spectrality of Bernoulli convolutions has been fully classified, the spectrality of the planar additive measure \(\mu_{r,s}\) cannot be directly deduced from the spectrality of its individual components. The construction of an exponential orthonormal basis in this two-dimensional setting requires not only spectra for the individual one-dimensional components, but also a precise alignment of their cross-component phases. On the other hand, it follows from a result of Lai, Liu and Prince \cite{LLP21} \footnote{For a general additive measure \(\eta = \frac12(\mu\times\delta_0+\delta_0\times\nu)\), if \(E(\Lambda)\) is a Fourier frame for \(L^2(\eta)\), then \(E(\pi_1(\Lambda))\) and \(E(\pi_2(\Lambda))\) are Fourier frames for \(L^2(\mu)\) and \(L^2(\nu)\), respectively, where \(\pi_1, \pi_2\) denote the coordinate projections onto the \(x\)- and \(y\)-axes (see Section 4 of~\cite{LLP21}).}  that the frame spectrality of \(\mu_{r,s}\) necessitates the frame spectrality of its one-dimensional component. Consequently, the failure of frame spectrality in either component immediately precludes the frame spectrality of $\mu_{r,s}$. In particular, recent works by de Dios Pont, Liehr and Taylor~\cite{DLT26} (and independently by Fu, Song and Wang~\cite{FSW26}) proved that the Bernoulli convolutions \(\mu_{1/N}\) for odd integers \(N \ge 3\) do not even admit Fourier frames, which is a stronger obstruction than the failure of spectrality. This motivates the following central question:

\begin{ques}
For which pairs \((r,s)\) is \(\mu_{r,s}\) a spectral measure? In particular,
what changes when \(r=1/(2k_1)\) and \(s=1/(2k_2)\) with \(k_1\ne k_2\)?
\end{ques}

Our first main result focuses on the symmetric case of equal even reciprocal contractions, i.e., \(r=s=1/(2k)\). For an even integer \(R=2k \ge 2\), our primary goal is to construct a spectrum of the planar measure \(\mu^{1/R}\). To achieve this, we first establish a new spectrum of the one-dimensional component \(\mu_{1/R}\), which then serves as a fundamental building block for the planar case.

We begin with the set \(\Gamma_R\), consisting of all finite sums of the form \(\sum_{j=0}^{N}\varepsilon_jR^j\) with \(\varepsilon_j\in\{0,1\}\). We then construct a spectrum of \(\mu_{1/R}\) as the disjoint union of a scaled copy of \(\Gamma_R\) and its translated reflection. Specifically, letting \(q_R=\frac{R(R-1)}4\), we define
\begin{equation}\label{eq:ST}
 S_R=(R-1)q_R\Gamma_R\ \cup\ (-q_R-(R-1)q_R\Gamma_R).
\end{equation}
It is worth emphasizing that this construction yields a new spectrum of \(\mu_{1/R}\), which cannot be deduced directly from any previously known constructions in the literature \cite{DFK24+,DHS09,DH16,FHW18,LW02}.

To construct a spectrum of the planar measure \(\mu^{1/R}\), we adjoin a translate of this one-dimensional spectrum \(S_R\) defined in \eqref{eq:ST} to itself, yielding
\[
 T_R=S_R\ \cup\ (S_R+(R-1)/4).
\]
The planar spectrum is then obtained by placing \(T_R\) on the anti-diagonal:
\[
 \Lambda_R=\{(t,-t):t\in T_R\}.
\]

With this explicit construction in hand, our first main theorem reads as follows.

\begin{thm}\label{thm:even}
For every \(k\ge1\), the one-dimensional measure \(\mu_{1/(2k)}\) admits a spectrum \(S_{2k}\), and the planar measure \(\mu^{1/(2k)}\) has a spectrum \(\Lambda_{2k}\).
\end{thm}

Combining Theorem \ref{thm:even} with the aforementioned established absence of Fourier frames for \(\mu^{1/R}\) when \(R\) is odd, we obtain the following parity dichotomy for the symmetric additive measures.

\begin{coro}\label{cor:parity-intro}
For every integer \(R\ge2\), $\mu^{1/R}$ is spectral if and only if $R$ is even.
\end{coro}

Naturally, one may ask what happens if the symmetry is broken. We show that the spectrality does not generally persist for unequal contractions, as demonstrated by the following counterexample.

\begin{thm}\label{thm:unequal}
The additive measure
\[
 \mu_{1/2,1/4}
 =\frac12(\mu_{1/2}\times\delta_0+\delta_0\times\mu_{1/4})
\]
is not a spectral measure.
\end{thm}

The paper is organized as follows. In Section~\ref{sec:prelim}, we introduce some basic definitions and preliminaries on spectral measures.  In Section~\ref{sec:S},  we prove Theorem~\ref{thm:even}. In Section~\ref{sec:unequal}, we prove
Theorem~\ref{thm:unequal}. In Section~\ref{sec:parity}, we prove Corollary \ref{cor:parity-intro}. In the last section, we give some examples and concluding remarks.

\section{Preliminaries}\label{sec:prelim}
Let $\sigma$ be a probability measure with compact support on $\R^d$. The Fourier transform of $\sigma$ is defined by
\[
 \hat\sigma(\xi)=\int e^{-2\pi i\langle\xi,x\rangle}\dd\sigma(x),
\]
It is easy to see that $E(\Lambda)$ is an orthogonal set in $L^2(\sigma)$ (in this case, we call $\Lambda$ an {\it orhtogonal set} of $\sigma$) if and only if
      \begin{equation}\label{eq:minisdwa}
      (\Lambda-\Lambda)\setminus\{0\}\subset \mathcal{Z}(\hat{\sigma}),
      \end{equation}
where $\mathcal{Z}(f)$ denotes the zero set of a function, i.e., $\mathcal{Z}(f)=\{\xi\in \R^d:f(\xi)=0\}$.

Let $R\ge2$ be an even integer. Equation
\eqref{eq:ifs} gives
\begin{equation}\label{eq:mask}
 \widehat{\mu_{1/R}}(\xi)=m(\xi/R) \widehat{\mu_{1/R}}(\xi/R),
\end{equation}
where  $m(\xi):=\tfrac12(1+e^{-4\pi i\xi})=e^{-2\pi i\xi}\cos(2\pi \xi)$. Since $\widehat{\mu_{1/R}}(0)=1$, iterating \eqref{eq:mask} to infinity yields the infinite product
\begin{equation}\label{eq:product}
\widehat{\mu_{1/R}}(\xi)=e^{-2\pi i\xi/(R-1)}\prod_{n=1}^{\infty}
                \cos\!\left(\frac{2\pi \xi}{R^n}\right).
\end{equation}
Therefore
\begin{equation}\label{eq:zeros}
 \mathcal{Z}(\widehat{\mu_{1/R}})=\left\{\frac{R^{k}}4(2j+1):k\geq1,\ j\in\Z\right\}.
\end{equation}

If $\eta=\frac12(\mu\times\delta_0+\delta_0\times\nu)$ for Borel probability
measures $\mu,\nu$ on $\R$, then
\begin{equation}\label{eq:additive-ft}
 \widehat{\eta}(\xi,\ell)=\frac12\bigl(\widehat\mu(\xi)+\widehat\nu(\ell)\bigr),
\end{equation}
as is immediate from the definition of the Fourier transform.

The following criterion is a universal test for a set $\Lambda\subset\R^d$, which is a sufficient and necessary condition to decide whether $\Lambda$ is an orthogonal set (or a spectrum) of $\sigma$.  Define
\begin{equation}\label{eq:JP}
 Q^\sigma_\Lambda(\xi)=\sum_{\lambda\in\Lambda}|\hat{\sigma}(\xi+\lambda)|^2,
\end{equation}
which plays a key role in the criterion.
\begin{thm}[{\cite{JP98}}]\label{thm:JP}
Let $\sigma$ be a Borel probability measure with compact support on $\mathbb R^d$ and let
$\Lambda\subset\mathbb R^d$ be a countable set. Then
\begin{enumerate}
    \item[\textup{(1)}] $\Lambda$ is an orthogonal set of $\sigma$ if and only if $Q^\sigma_\Lambda(\xi)\le1$ for all $\xi\in\mathbb R^d$;
    \item[\textup{(2)}] $\Lambda$ is a spectrum of $\sigma$ if and only if $Q^\sigma_\Lambda(\xi)=1$ for all $\xi\in\mathbb R^d$;
    \item[\textup{(3)}] $Q^\sigma_\Lambda(\xi)$ has an entire analytic extension to $\mathbb C^n$ if $\Lambda$ is an orthogonal set of $\sigma$.
\end{enumerate}
\end{thm}

\section{Proof of Theorem \ref{thm:even}}\label{sec:S}
In this section, we prove Theorem \ref{thm:even}. We first prove the first assertion of Theorem \ref{thm:even}, i.e., $S_R$ defined in \eqref{eq:ST} is a spectrum of $\mu_{1/R}$. By Theorem \ref{thm:JP}, this reduces to proving that $Q_R(\xi):=Q^{\mu_{1/R}}_{S_R}(\xi)\equiv1$ for all $\xi\in\R$. The proof consists of the following three steps:

\begin{enumerate}
  \item[(i)] The self-similar recursion
  \[
    S_R=RS_R\cup\bigl(RS_R+(R-1)q_R\bigr)
  \]
  is established in Lemma~\ref{lem:rec}; together with \eqref{eq:zeros}, it yields the orthogonality of \(\{e_s:s\in S_R\}\) in Lemma~\ref{lem:orth}.

  \item[(ii)] Using this recursion and the mask relation \eqref{eq:mask},we derive in \eqref{eq:transfer} the two-branch functional
  equation for the spectral function
  \[
    Q_R(t)=\sum_{s\in S_R}
    \bigl|\widehat{\mu_{1/R}}(t+s)\bigr|^2,
  \]
  namely
  \[
    Q_R(t)
    =\cos^2\!\Bigl(\frac{2\pi t}{R}\Bigr)Q_R\!\Bigl(\frac tR\Bigr)
    +\sin^2\!\Bigl(\frac{2\pi t}{R}\Bigr)
     Q_R\!\Bigl(\frac{t+\ell_R}{R}\Bigr).
  \]

  \item[(iii)] The maximum-principle argument in Proposition~\ref{prop:complete} then gives \(Q_R\equiv1\), so \(S_R\) is a spectrum of \(\mu_{1/R}\) by Theorem~\ref{thm:JP}.
\end{enumerate}

To simplify notation define $\ell_R=\frac{R(R-1)^2}{4}$.
\begin{lemma}\label{lem:rec}
Let $R\ge 2$ be an even integer.
Then the set $S_R$ satisfies
\begin{equation}\label{eq:recursion}
S_R = R S_R \cup \big(R S_R + \ell_R\big),
\end{equation}
with $0,-q_R\in S_R$, where the union is disjoint.
\end{lemma}
\begin{proof}
From the definition of $\Gamma_R$ it is
immediate that
\[
 \Gamma_R=R\Gamma_R\,\cup\,(R\Gamma_R+1).
\]
Applying
this to the two components of $S_R$ in \eqref{eq:ST} gives
\begin{align*}
 \ell_R\Gamma_R
 &=R\ell_R\Gamma_R\ \cup\ (R\ell_R\Gamma_R+\ell_R),\\
 -q_R-\ell_R\Gamma_R
 &=(-q_R-R\ell_R\Gamma_R)\ \cup\ (-Rq_R-R\ell_R\Gamma_R),
\end{align*}
where the second line uses $-q_R-\ell_R=-Rq_R$. Since $-Rq_R+\ell_R=-q_R$, the four displayed pieces regroup as
\begin{align}\label{eq:ellaer}
 S_R &= \bigl(R\ell_R\Gamma_R\cup(-Rq_R-R\ell_R\Gamma_R)\bigr)\cup\ \bigl((R\ell_R\Gamma_R+\ell_R)\cup(-q_R-R\ell_R\Gamma_R)\bigr) \nonumber \\
     &= RS_R\ \cup\ (RS_R+\ell_R),
\end{align}
which is \eqref{eq:recursion}.

To see the above two branches in \eqref{eq:ellaer} are disjoint, note that $\Gamma_R\subset\mathbb Z$ and
\[
 \ell_R=\frac{R(R-1)^2}{4}\in\frac12\mathbb Z,
 \qquad
 q_R=\frac{R(R-1)}4\in\frac12\mathbb Z,
\]
we have $S_R\subset\frac12\mathbb Z$. If $Rs=Rt+(R-1)q_R$ for some
$s,t\in S_R$, then
\[
 s-t=\frac{(R-1)q_R}{R}=\frac{\ell_R}{R}
      =\frac{(R-1)^2}{4}.
\]
As $R$ is even, $R-1$ is odd, hence $(R-1)^2\equiv1\pmod4$ and therefore
\[
 \frac{(R-1)^2}{4}\in\mathbb Z+\frac14.
\]
This is not an element of $\frac12\mathbb Z$, contradicting
$s-t\in\frac12\mathbb Z$. Thus the two branches are disjoint.

Finally, $0$ and $-q_R$ are the zero-digit members of the two components
of \eqref{eq:ST}: namely,
\[
 0=\ell_R\cdot0\in\ell_R\Gamma_R,
 \qquad
 -q_R=-q_R-\ell_R\cdot0\in -q_R-\ell_R\Gamma_R.
\]
Hence $0,-q_R\in S_R$, as asserted.
\end{proof}

\begin{lemma}\label{lem:orth}
The family $\{e_s:s\in S_R\}$ is orthonormal in $L^2(\mu_{1/R})$.
\end{lemma}
\begin{proof}
By \eqref{eq:minisdwa}, it suffices to prove
\[
 (S_R-S_R)\setminus\{0\}\subset
 \mathcal Z\bigl(\widehat{\mu_{1/R}}\bigr).
\]
For any  $s,t$ in $S_R$ with $s\neq t$. By Lemma~\ref{lem:rec}, write
\(s=Rs_1+a\) and \(t=Rt_1+b\), where $a,b\in\{0,\ell_R\}$ and $s_1,t_1\in S_R$.

If $a=b$, then $s-t=R(s_1-t_1)$ with $s_1\neq t_1$. Replacing $(s,t)$ by
$(s_1,t_1)$ and iterating, we eventually reach the case $a\neq b$. Indeed,
otherwise $s-t=R^n(s_n-t_n)$ for every $n$, with $s_n\neq t_n$ and
$s_n,t_n\in S_R\subset\frac12\mathbb Z$ (as shown in the proof of
Lemma~\ref{lem:rec}). Hence $|s_n-t_n|\ge\frac12$, which would force
\[
 |s-t|\ge \frac{R^n}{2}\to\infty,
\]
a contradiction.

Thus we may assume $a\neq b$, so that $\{a,b\}=\{0,\ell_R\}$ and hence
$a-b=\pm\ell_R$. Since $s_1-t_1\in\frac12\mathbb Z$, there exists
$j\in\mathbb Z$ such that $s_1-t_1=j/2$. Therefore
\[
 s-t=R(s_1-t_1)+(a-b)
     =R\cdot\frac j2\pm\ell_R
     =\frac R4\bigl(2j\pm(R-1)^2\bigr).
\]
The integer $2j\pm(R-1)^2$ is odd, so
$s-t\in\mathcal Z(\widehat{\mu_{1/R}})$ by \eqref{eq:zeros}. Hence, we complete the proof.
\end{proof}

By Theorem~\ref{thm:JP}, it remains to prove \(Q_R\equiv1\), where
\begin{equation}\label{eq:Qdef}
 Q_R(t):=\sum_{s\in S_R}
 \bigl|\widehat{\mu_{1/R}}(t+s)\bigr|^2.
\end{equation}
Since \(S_R\) is orthonormal by Lemma~\ref{lem:orth}, Theorem~\ref{thm:JP}(1) gives
\begin{equation}\label{eq:Qle1}
 0\le Q_R(t)\le1\qquad(t\in\R).
\end{equation}
By Lemma~\ref{lem:rec}, \(0,-q_R\in S_R\). Hence at \(t=0\) the summand \(s=0\), and at \(t=q_R\) the summand \(s=-q_R\), contribute
\[
 \bigl|\widehat{\mu_{1/R}}(0)\bigr|^2=1.
\]
Since all summands are nonnegative and \(Q_R\le1\), this forces
\begin{equation}\label{eq:anchor}
 Q_R(0)=Q_R(q_R)=1.
\end{equation}

The key ingredient is a self-similarity identity for \(Q_R\). Substituting the recursion \eqref{eq:recursion} into \eqref{eq:Qdef} and applying the mask relation \eqref{eq:mask}, we obtain
\begin{equation}\label{eq:transfer}
 Q_R(t)=
 \cos^2\!\Bigl(\frac{2\pi t}{R}\Bigr)Q_R\!\Bigl(\frac tR\Bigr)
 +\sin^2\!\Bigl(\frac{2\pi t}{R}\Bigr)
 Q_R\!\Bigl(\frac{t+\ell_R}{R}\Bigr).
\end{equation}
Indeed, for \(s=Rs'\) we have
\[
 \widehat{\mu_{1/R}}(t+Rs')
 =m(t/R+s')\widehat{\mu_{1/R}}(t/R+s').
\]
Since \(|m|^2=\cos^2(2\pi\,\cdot)\) is \(\frac12\)-periodic and \(S_R\subset\frac12\mathbb Z\), we get
\[
 |m(t/R+s')|^2=\cos^2\!\Bigl(\frac{2\pi t}{R}\Bigr).
\]
For \(s=Rs'+\ell_R\), the argument of \(m\) is \(t/R+s'+\ell_R/R\). Since
\[
 \frac{\ell_R}{R}=\frac{(R-1)^2}{4}\in\mathbb Z+\frac14
\]
(as \(R\) is even) and \(S_R\subset\frac12\mathbb Z\), this turns \(\cos^2\) into \(\sin^2\), while \(\widehat{\mu_{1/R}}\) is evaluated at \((t+\ell_R)/R+s'\). Summing over \(s'\in S_R\) gives \eqref{eq:transfer}.

Write
\begin{equation}\label{eq:sigma}
 \sigma_0(t)=\frac tR,\qquad \sigma_1(t)=\frac{t+\ell_R}{R},
 \qquad
 w_0(t)=\cos^2\!\Bigl(\frac{2\pi t}{R}\Bigr),\qquad
 w_1(t)=\sin^2\!\Bigl(\frac{2\pi t}{R}\Bigr).
\end{equation}
Then \(w_0,w_1\ge0\) and \(w_0+w_1=1\). Consequently, subtracting \eqref{eq:transfer} from \(1\), the nonnegative function \(1-Q_R\) satisfies the same convex identity:
\begin{equation}\label{eq:transferh}
 1-Q_R(t)=w_0(t)\bigl(1-Q_R(\sigma_0(t))\bigr)
 +w_1(t)\bigl(1-Q_R(\sigma_1(t))\bigr).
\end{equation}

\begin{proposition}\label{prop:complete}
For every even $R\ge2$, \(Q_R(t)=1\) for all \(t\in\R\); hence \(S_R\) is a spectrum of \(\mu_{1/R}\).
\end{proposition}

\begin{proof}
We prove \(Q_R\equiv1\). Since \(S_R\) is orthogonal, Theorem~\ref{thm:JP}(3) shows that \(Q_R\), and hence \(1-Q_R\), has an entire analytic extension; in particular \(1-Q_R\) is continuous.

Put \(I_R=[0,q_R]\). Both maps in \eqref{eq:sigma} preserve \(I_R\):
\[
 \sigma_0(I_R)=\Bigl[0,\frac{q_R}{R}\Bigr]\subset I_R,\qquad
 \sigma_1(0)=\frac{\ell_R}{R}\in[0,q_R],\qquad
 \sigma_1(q_R)=q_R,
\]
where \(\sigma_1(q_R)=q_R\) uses \(\ell_R=(R-1)q_R\).

Let
\[
 M_*=\max_{t\in I_R}\bigl(1-Q_R(t)\bigr),
\]
and suppose \(M_*>0\). Choose \(t_0\in I_R\) such that
\[
 1-Q_R(t_0)=M_*.
\]
By \eqref{eq:anchor}, \(Q_R(0)=Q_R(q_R)=1\), so \(t_0\in(0,q_R)\).

\emph{Case 1: \(w_0(t_0)>0\).}
By \eqref{eq:transferh},
\[
 1-Q_R(t_0)
 =w_0(t_0)\bigl(1-Q_R(\sigma_0(t_0))\bigr)
 +w_1(t_0)\bigl(1-Q_R(\sigma_1(t_0))\bigr).
\]
The left-hand side equals \(M_*\), while both terms in parentheses are at most \(M_*\). Since \(w_0(t_0)>0\), we must have
\[
 1-Q_R(\sigma_0(t_0))=M_*,
\]
i.e.\ \(1-Q_R(t_0/R)=M_*\). Now
\[
 0<\frac{t_0}{R}<\frac{q_R}{R}=\frac{R-1}{4}<\frac R4,
\]
and \(R/4\) is the smallest positive zero of \(w_0\), so \(w_0(t_0/R^n)>0\) for every \(n\ge1\). Iterating the same argument gives
\[
 1-Q_R(t_0/R^n)=M_*\qquad(n\ge1).
\]
Letting \(n\to\infty\) and using continuity at \(0\), we obtain
\[
 1-Q_R(0)=M_*>0,
\]
contradicting \(Q_R(0)=1\).

\emph{Case 2: \(w_0(t_0)=0\).}
Then \(w_1(t_0)=1\) and
\begin{equation}\label{eq:exceptional}
 t_0=\frac R4 j,\qquad j\ \text{odd},\quad 1\le j\le R-3
\end{equation}
(this range is empty for \(R=2\)). Identity \eqref{eq:transferh} forces
\[
 1-Q_R(\sigma_1(t_0))=M_*,
\]
and
\begin{equation}\label{eq:exceptionalmapped}
 \sigma_1(t_0)=\frac{j+(R-1)^2}{4}.
\end{equation}
We claim \(w_0(\sigma_1(t_0))>0\). If \(w_0(\sigma_1(t_0))=0\), then
\(\sigma_1(t_0)=Rj'/4\) for some odd integer \(j'\), so
\[
 j+(R-1)^2=Rj'.
\]
Reducing modulo \(R\) gives \(j+1\equiv0\pmod R\), impossible since
\(2\le j+1\le R-2\). Thus \(w_0(\sigma_1(t_0))>0\), and Case~1 applied to
\(\sigma_1(t_0)\in I_R\) gives a contradiction.

Both cases are impossible, so \(M_*=0\), and therefore
\begin{equation}\label{eq:Q-on-I}
 Q_R(t)=1\qquad(t\in I_R).
\end{equation}
Thus the entire function \(1-Q_R\) vanishes on the interval \(I_R\). By the identity theorem, \(1-Q_R\equiv0\) on \(\R\), i.e.\ \(Q_R\equiv1\). Theorem~\ref{thm:JP} then shows that \(S_R\) is a spectrum of \(\mu_{1/R}\).
\end{proof}

\begin{remark}[Why the parameter $\ell_R$ is not arbitrary]
The three identities
\[
 \ell_R=(R-1)q_R,\qquad \ell_R/R\in\Z+\tfrac14,\qquad S_R\subset D_R\Z
\]
perform three different jobs: they give the exact two-branch recursion, the
complementary mask weights, and a lattice compatible with reflection. Taking the
canonical one-dimensional spectrum from Dai's theorem without these restrictions
would not prove the planar assertion of the next section.
\end{remark}

We now prove the second assertion of Theorem~\ref{thm:even}. Recall that
\[
 T_R=S_R\cup\left(S_R+\frac{R-1}{4}\right),
 \qquad
 \Lambda_R=\{(t,-t):t\in T_R\}.
\]
We compute the Jorgensen--Pedersen spectral function of \(\Lambda_R\) and show that it equals \(1\) everywhere; by Theorem~\ref{thm:JP} this proves that \(\Lambda_R\) is a spectrum of \(\mu^{1/R}\).

Since
\begin{equation}\label{eq:rhohat}
 \widehat{\mu^{1/R}}(\xi,\eta)
 =\tfrac12\bigl(\widehat{\mu_{1/R}}(\xi)+\widehat{\mu_{1/R}}(\eta)\bigr),
\end{equation}
we have
\begin{equation}\label{eq:Qr}
 Q_{\Lambda_R}^{\mu^{1/R}}(\xi,\eta)
 =\frac14\sum_{t\in T_R}
 \bigl|\widehat{\mu_{1/R}}(\xi+t)+\widehat{\mu_{1/R}}(\eta-t)\bigr|^2.
\end{equation}
By Proposition~\ref{prop:complete}, \(Q_R\equiv1\). Hence the two projection sums are
\begin{equation}\label{eq:twotight}
 \sum_{t\in T_R}\bigl|\widehat{\mu_{1/R}}(\xi+t)\bigr|^2
 =Q_R(\xi)+Q_R\left(\xi+\frac{R-1}{4}\right)=2,
\end{equation}
and similarly
\begin{equation}\label{eq:twotight2}
 \sum_{t\in T_R}\bigl|\widehat{\mu_{1/R}}(\eta-t)\bigr|^2=2.
\end{equation}

It remains to see that the cross terms cancel. We first record the Parseval
identity underlying the computation. Since \(\{e_s:s\in S_R\}\) is an
orthonormal basis of \(L^2(\mu_{1/R})\), for every \(\xi,\eta\in\R\),
\begin{equation}\label{eq:parseval}
 \sum_{s\in S_R}
 \widehat{\mu_{1/R}}(\xi+s)\,\widehat{\mu_{1/R}}(\eta-s)
 =\widehat{\mu_{1/R}}(\xi+\eta).
\end{equation}
Indeed, \(\widehat{\mu_{1/R}}(\xi+s)=\ip{e_{-\xi}}{e_s}\) and
\(\widehat{\mu_{1/R}}(\eta-s)=\ip{e_s}{e_\eta}\), so by Parseval
\[
 \sum_{s\in S_R}\widehat{\mu_{1/R}}(\xi+s)\,\widehat{\mu_{1/R}}(\eta-s)
 =\sum_{s\in S_R}\ip{e_{-\xi}}{e_s}\ip{e_s}{e_\eta}
 =\ip{e_{-\xi}}{e_\eta}
 =\widehat{\mu_{1/R}}(\xi+\eta).
\]

Next, write
\[
 \Phi(u)=\prod_{n=1}^{\infty}\cos\!\Bigl(\frac{2\pi u}{R^n}\Bigr),
\]
so that \(\Phi\) is real and even and, by \eqref{eq:product},
\begin{equation}\label{eq:phase}
 \widehat{\mu_{1/R}}(u)=e^{-2\pi iu/(R-1)}\Phi(u).
\end{equation}
We also use that \(S_R\subset\frac{R-1}{2}\Z\). Indeed
\(\ell_R=\frac{R-1}{2}\cdot\frac{R(R-1)}2\) and
\(q_R=\frac{R-1}{2}\cdot\frac R2\), with
\(\frac{R(R-1)}2,\frac R2\in\Z\) because \(R\) is even; hence
\(S_R=\ell_R\Gamma_R\cup(-q_R-\ell_R\Gamma_R)\subset\frac{R-1}{2}\Z\).

Define
\begin{equation}\label{eq:K0}
 K_0(\xi,\eta)=\sum_{s\in S_R}
 \widehat{\mu_{1/R}}(\xi+s)
 \overline{\widehat{\mu_{1/R}}(\eta-s)}.
\end{equation}
The sum converges absolutely by Cauchy--Schwarz and Parseval:
\[
 \sum_{s\in S_R}
 \bigl|\widehat{\mu_{1/R}}(\xi+s)
 \widehat{\mu_{1/R}}(\eta-s)\bigr|\le1.
\]
By \eqref{eq:phase}, since \(\Phi\) is real and
\(e^{-4\pi is/(R-1)}=1\) for \(s\in\frac{R-1}{2}\Z\), we obtain
\[
 K_0(\xi,\eta)
 =e^{-2\pi i(\xi-\eta)/(R-1)}
   \sum_{s\in S_R}\Phi(\xi+s)\,\Phi(\eta-s).
\]
Inserting \eqref{eq:phase} into \eqref{eq:parseval} and cancelling the common
phase factor \(e^{-2\pi i(\xi+\eta)/(R-1)}\) yields
\[
 \sum_{s\in S_R}\Phi(\xi+s)\,\Phi(\eta-s)=\Phi(\xi+\eta).
\]
Therefore
\begin{equation}\label{eq:Kvalue}
 K_0(\xi,\eta)
 =e^{-2\pi i(\xi-\eta)/(R-1)}\Phi(\xi+\eta)
 =e^{4\pi i\eta/(R-1)}\widehat{\mu_{1/R}}(\xi+\eta).
\end{equation}
The shifted branch gives
\begin{align}\label{eq:K1}
 K_1(\xi,\eta)
 &=\sum_{s\in S_R}
 \widehat{\mu_{1/R}}\left(\xi+s+\frac{R-1}{4}\right)
 \overline{\widehat{\mu_{1/R}}\left(\eta-s-\frac{R-1}{4}\right)}
 \notag\\
 &=K_0\left(\xi+\frac{R-1}{4},\eta-\frac{R-1}{4}\right)\notag\\
 &=e^{4\pi i(\eta-(R-1)/4)/(R-1)}\widehat{\mu_{1/R}}(\xi+\eta)
 =-K_0(\xi,\eta),
\end{align}
because
\[
 e^{-4\pi i(R-1)/(4(R-1))}=e^{-\pi i}=-1.
\]
Therefore the cross terms cancel exactly:
\begin{equation}\label{eq:Kcancel}
 \sum_{t\in T_R}
 \widehat{\mu_{1/R}}(\xi+t)
 \overline{\widehat{\mu_{1/R}}(\eta-t)}
 =K_0(\xi,\eta)+K_1(\xi,\eta)=0.
\end{equation}
Combining \eqref{eq:Qr}, \eqref{eq:twotight}, \eqref{eq:twotight2}, and \eqref{eq:Kcancel} gives
\begin{equation}\label{eq:Qmain}
 Q_{\Lambda_R}^{\mu^{1/R}}(\xi,\eta)
 =\frac14\left(2+2+2\operatorname{Re}0\right)=1
 \qquad((\xi,\eta)\in\R^2).
\end{equation}
By Theorem~\ref{thm:JP}, \(\Lambda_R\) is a spectrum of \(\mu^{1/R}\). This proves Theorem~\ref{thm:even}.

\section{Proof of Theorem \ref{thm:unequal}}\label{sec:unequal}

In this section we prove Theorem~\ref{thm:unequal}, namely that the additive measure \(\mu_{1/2,1/4}\) is not spectral. The two components have different analytic types: \(\mu_{1/2}\) is the normalized Lebesgue measure on \([0,2]\), whereas \(\mu_{1/4}\) is singular continuous. Crucially, the argument below imposes no restriction on the geometry of a putative spectrum \(\Lambda\); in particular, it never assumes that \(\Lambda\) lies on a line or splits into finitely many branches.

We use the Fourier identities
\begin{align}
 \widehat{\mu_{1/2}}(a)
 &=e^{-2\pi ia}\frac{\sin(2\pi a)}{2\pi a},\label{eq:F2}\\
 \widehat{\mu_{1/4}}(b)
 &=\prod_{j=1}^{\infty}
   e^{-2\pi ib/4^j}\cos(2\pi b/4^j).\label{eq:F4}
\end{align}
The quotient at \(a=0\) is understood as \(1\). As before, the additive measure has Fourier transform
\begin{equation}\label{eq:rho2}
 \widehat{\mu_{1/2,1/4}}(\xi,\eta)
 =\tfrac12\bigl(
 \widehat{\mu_{1/2}}(\xi)+\widehat{\mu_{1/4}}(\eta)
 \bigr).
\end{equation}
Our main analytic input is a sharp bound on the part of the range of \(\widehat{\mu_{1/4}}\) that lies in the closed left half-plane. Its proof is self-contained and does not rely on numerical evidence or on any classification of the frequency set.

\begin{proposition}\label{prop:range}
For every \(b\in\R\),
\begin{equation}\label{eq:range}
 \Ree\widehat{\mu_{1/4}}(b)\leq0
 \quad\Longrightarrow\quad
 \bigl|\widehat{\mu_{1/4}}(b)\bigr|<\frac2\pi.
\end{equation}
\end{proposition}

To prove Proposition~\ref{prop:range}, we need the following two elementary lemmas.

\begin{lemma}\label{lem:residues}
For each \(j\geq1\), choose an integer \(N_j\) such that
\[
 r_j=\frac{2b}{4^j}-N_j\in\left[-\tfrac12,\tfrac12\right].
\]
Then \(N_j=0\) for all sufficiently large \(j\), and
\begin{equation}\label{eq:residueproduct}
 \widehat{\mu_{1/4}}(b)
 =e^{-\pi i\sum_{j\geq1}r_j}\prod_{j\geq1}\cos(\pi r_j).
\end{equation}
The sum converges absolutely. Moreover,
\begin{equation}\label{eq:carry}
 k_j:=4r_{j+1}-r_j=N_j-4N_{j+1}\in\Z,
\end{equation}
and only finitely many \(k_j\) are nonzero.
\end{lemma}

\begin{proof}
Since \(2b/4^j\to0\), its nearest integer is eventually zero. For any integer \(N\),
\[
 e^{-\pi i(N+r)}\cos(\pi(N+r))
 =e^{-\pi ir}\cos(\pi r).
\]
Applying this to each factor of \eqref{eq:F4}, we obtain \eqref{eq:residueproduct}.
Eventually \(r_j=2b/4^j\), so the tail is geometric and all asserted limits
and identities follow. At a nearest-integer tie one factor is zero; such points
cannot satisfy \(|\widehat{\mu_{1/4}}(b)|\geq2/\pi\).
\end{proof}

\begin{lemma}\label{lem:onecarry}
If \(|\widehat{\mu_{1/4}}(b)|\geq2/\pi\), then \(|r_j|<1/3\) for every \(j\)
and the sequence \((k_j)\) is either identically zero or has exactly one nonzero
entry, which is \(+1\) or \(-1\).
\end{lemma}

\begin{proof}
Every factor in the nonnegative product \eqref{eq:residueproduct} is at most one.
Consequently
\[
 \cos(\pi r_j)\geq |\widehat{\mu_{1/4}}(b)|\geq2/\pi>1/2,
\]
which gives \(|r_j|<1/3\). Thus \(|k_j|<5/3\), so \(k_j\in\{-1,0,1\}\).

For \(|u|<\pi/2\), the inequality \(\log\cos u\leq-u^2/2\) follows by integrating
\(\tan u\geq u\) on \([0,|u|]\). Hence
\begin{equation}\label{eq:gaussianbound}
 |\widehat{\mu_{1/4}}(b)|
 \leq\exp\left(-\frac{\pi^2}{2}\sum_{j\geq1}r_j^2\right).
\end{equation}
Suppose that two carries are nonzero, at indices \(j<\ell\).
If \(\ell\geq j+2\), Cauchy--Schwarz gives
\[
 1=(4r_{j+1}-r_j)^2\leq17(r_j^2+r_{j+1}^2),
\]
and the same estimate for \(\ell\), involving disjoint variables. Therefore
\(\sum r_i^2\geq2/17\).

If \(\ell=j+1\), put \(v=(r_j,r_{j+1},r_{j+2})^T\) and
\[
 A=\begin{pmatrix}-1&4&0\\0&-1&4\end{pmatrix}.
\]
The eigenvalues of \(AA^T=\begin{pmatrix}17&-4\\-4&17\end{pmatrix}\) are \(13\)
and \(21\). Since \(Av=(k_j,k_{j+1})^T\) has squared norm two,
\[
 2\leq21\|v\|^2,
 \qquad
 \sum_{i\geq1}r_i^2\geq\frac2{21}.
\]
In both cases \eqref{eq:gaussianbound} yields
\begin{equation}\label{eq:constant}
 |\widehat{\mu_{1/4}}(b)|\leq e^{-\pi^2/21}<2/\pi,
\end{equation}
a contradiction. The last strict inequality is elementary and does not rely on
decimal evaluation: \(\pi^2>49/5\) and \(\pi<22/7\) imply
\[
 e^{\pi^2/21}>1+\frac7{15}+\frac12\left(\frac7{15}\right)^2
 =\frac{709}{450}>\frac{11}{7}>\frac\pi2.
\]
Thus there is at most one nonzero carry.
\end{proof}

\begin{proof}[Proof of Proposition~\ref{prop:range}]
Suppose for a contradiction that
\[
 \Ree\widehat{\mu_{1/4}}(b)\leq0
 \quad\text{and}\quad
 \bigl|\widehat{\mu_{1/4}}(b)\bigr|\geq\frac2\pi.
\]
By Lemma~\ref{lem:onecarry}, put \(r=r_1\).

\emph{No nonzero carry.}
In this case \(r_j=r/4^{j-1}\), so
\[
 \sum_{j\geq1}r_j=\frac{4r}{3}\in\left(-\frac49,\frac49\right).
\]
Since the product of cosines in \eqref{eq:residueproduct} is strictly positive,
\(\widehat{\mu_{1/4}}(b)\) has positive real part, a contradiction.

\emph{One nonzero carry.}
Replacing \(b\) by \(-b\) conjugates \(\widehat{\mu_{1/4}}(b)\) and changes all
residue signs. We may therefore assume that the carry is \(k_m=1\) at some
\(m\geq1\). The recurrence \eqref{eq:carry} gives
\begin{align}
 r_j&=r/4^{j-1} &&(1\leq j\leq m),\nonumber\\
 r_{m+1}&=\tfrac14+r/4^m,\label{eq:singlecarry}\\
 r_{m+1+q}&=4^{-q}r_{m+1} &&(q\geq0).\nonumber
\end{align}
Summing,
\[
 \sum_{j\geq1}r_j=\frac{4r+1}{3}\in\left(-\frac19,\frac79\right).
\]
On this interval, \(\cos(\pi s)\leq0\) implies \(s\geq1/2\). Thus
\(\Ree\widehat{\mu_{1/4}}(b)\leq0\) forces
\begin{equation}\label{eq:rthreshold}
 r\geq\tfrac18.
\end{equation}

Define the positive product
\begin{equation}\label{eq:C}
 C(x)=\prod_{q=0}^{\infty}\cos\left(\frac{\pi x}{4^q}\right)
 \qquad\left(0\leq x<\frac12\right).
\end{equation}
It is nonincreasing on this interval. Put \(h=r/4^m>0\).
By \eqref{eq:singlecarry},
\begin{equation}\label{eq:quotient}
 \bigl|\widehat{\mu_{1/4}}(b)\bigr|
 =\frac{C(r)}{C(h)}C\left(\tfrac14+h\right).
\end{equation}
Here \(r<1/3\), \(m\geq1\), and hence \(0<h<1/12\).
For every \(q\geq0\), both angles below and their sum are less than \(\pi/2\),
and
\[
 \cos\left(\frac{\pi}{4^{q+1}}+\frac{\pi h}{4^q}\right)
 <\cos\left(\frac{\pi}{4^{q+1}}\right)
   \cos\left(\frac{\pi h}{4^q}\right).
\]
Multiplying gives
\begin{equation}\label{eq:strictproduct}
 C\left(\tfrac14+h\right)<C\left(\tfrac14\right)C(h).
\end{equation}
Strictness persists for the infinite products: the ratio of the first factors
is strictly less than one, and every later ratio is at most one; all three
products are positive. Equations \eqref{eq:rthreshold}--\eqref{eq:strictproduct}
imply
\[
 \bigl|\widehat{\mu_{1/4}}(b)\bigr|
 <C(r)C\left(\tfrac14\right)
 \leq C\left(\tfrac18\right)C\left(\tfrac14\right).
\]
Finally, interlacing the two products and using repeated angle doubling,
\begin{equation}\label{eq:sincproduct}
 C\left(\tfrac18\right)C\left(\tfrac14\right)
 =\prod_{n=2}^{\infty}\cos\left(\frac\pi{2^n}\right)
 =\frac{\sin(\pi/2)}{\pi/2}
 =\frac2\pi.
\end{equation}
This contradicts the assumed lower bound and proves the proposition.
\end{proof}

\begin{remark}[The constant is sharp as a supremum]\label{rem:sharp}
Let \(b_M=2\cdot4^M+1/4\) for \(M\geq0\). Its residues have \(r_1=1/8\)
and a single carry at \(m=M+1\). Therefore \(\widehat{\mu_{1/4}}(b_M)\) is
negative imaginary, and with \(h_M=1/(8\cdot4^{M+1})\),
\[
 \bigl|\widehat{\mu_{1/4}}(b_M)\bigr|
 =\frac{C(1/8)}{C(h_M)}C(1/4+h_M)<\frac2\pi,
 \qquad
 \bigl|\widehat{\mu_{1/4}}(b_M)\bigr|\to\frac2\pi.
\]
Thus the strict bound at every finite frequency is essential; replacing the
constant by a uniformly smaller one would be false.
\end{remark}

A first consequence of Proposition~\ref{prop:range} is a zero-free vertical
strip for the Fourier transform of \(\mu_{1/2,1/4}\).

\begin{corollary}\label{cor:strip}
For all \(a,b\in\R\),
\begin{equation}\label{eq:strip}
 |a|\leq\tfrac14
 \quad\Longrightarrow\quad
 \widehat{\mu_{1/2}}(a)+\widehat{\mu_{1/4}}(b)\neq0.
\end{equation}
Consequently, any orthogonal exponential set for \(\mu_{1/2,1/4}\) satisfies
\begin{equation}\label{eq:separation}
 |a-a'|>\tfrac14
 \quad\text{for distinct }(a,b),(a',b')\in\Lambda.
\end{equation}
\end{corollary}

\begin{proof}
For \(|a|\leq1/4\), formula \eqref{eq:F2} shows that
\[
 \Ree\widehat{\mu_{1/2}}(a)\geq0,
 \qquad
 \bigl|\widehat{\mu_{1/2}}(a)\bigr|\geq\frac2\pi.
\]
For example, the latter inequality follows from the chord bound
\(\sin u\geq2u/\pi\) on \([0,\pi/2]\). If
\(\widehat{\mu_{1/4}}(b)=-\widehat{\mu_{1/2}}(a)\), its real part is
nonpositive and its modulus is at least \(2/\pi\), contradicting
Proposition~\ref{prop:range}. Apply this to the difference of two orthogonal
frequencies, using \eqref{eq:rho2}, to obtain \eqref{eq:separation}.
In particular, horizontal projection is injective.
\end{proof}

We shall also need the following elementary packing lemma for tight frames on
an interval.

\begin{lemma}\label{lem:packing}
Let \(\mu_L=L^{-1}\mathbf1_{[0,L]}\,\mathrm dx\), where \(L>0\).
Suppose the exponentials with frequencies in \(A\subset\R\) form a tight frame
of bound \(c>0\) for \(L^2(\mu_L)\), with each exponential having norm one.
If distinct points of \(A\) have separation at least \(\delta=1/(cL)\),
then \(A=a_0+\delta\Z\) for some \(a_0\in\R\).
\end{lemma}

\begin{proof}
Let \(K(t)=|\widehat{\mu_L}(t)|^2\). Applying the tight-frame identity to an
exponential gives
\begin{equation}\label{eq:Kidentity}
 \sum_{a\in A}K(t-a)=c\qquad(t\in\R).
\end{equation}
The sinc formula, or Plancherel's theorem, gives
\(\int_\R K(t)\dd t=1/L\). Also, \(K\geq0\) and \(K>0\) except at a discrete set.

Put \(I=[-\delta/2,\delta/2]\) and \(h=\sum_{a\in A}\mathbf1_{a+I}\).
The spacing hypothesis implies \(0\leq h\leq1\) almost everywhere.
Tonelli's theorem and \eqref{eq:Kidentity} give, for each \(t\),
\[
 (h*K)(t)
 =\int_I\sum_{a\in A}K(t-u-a)\dd u
 =c\delta=\frac1L=(1*K)(t).
\]
Therefore \((1-h)*K=0\). Both factors are nonnegative, and \(K\) is positive
almost everywhere. Evaluating this convolution at any one \(t\) forces
\(1-h=0\) almost everywhere.

It follows that the length-\(\delta\) intervals centered at \(A\) cover the line
almost everywhere and have disjoint interiors. The uniform spacing makes \(A\)
locally finite; covering forces it to be unbounded in both directions.
Consecutive centers must have distance exactly \(\delta\), since a larger gap
would leave an open interval uncovered. Thus \(A\) is the claimed lattice.
\end{proof}

\begin{proof}[Proof of Theorem~\ref{thm:unequal}]
Assume that \(\Lambda\) is a spectrum of \(\mu_{1/2,1/4}\). Since neither
component has an atom at zero, the two coordinate axes intersect in a
\(\mu_{1/2,1/4}\)-null set, and a function
\(h\in L^2(\mu_{1/2,1/4})\) is determined by its restrictions
\(f=h(\cdot,0)\in L^2(\mu_{1/2})\) and \(g=h(0,\cdot)\in L^2(\mu_{1/4})\)
to the two axes, with norm
\[
 \|h\|_{L^2(\mu_{1/2,1/4})}^2
 =\tfrac12\|f\|_{L^2(\mu_{1/2})}^2
 +\tfrac12\|g\|_{L^2(\mu_{1/4})}^2.
\]
For \(f\in L^2(\mu_{1/2})\), Parseval applied to \((f,0)\) gives
\[
 \frac12\|f\|_{L^2(\mu_{1/2})}^2
 =\sum_{(a,b)\in\Lambda}
 \left|\frac12\langle f,e_a\rangle_{L^2(\mu_{1/2})}\right|^2.
\]
By Corollary~\ref{cor:strip}, horizontal projection is injective. Thus
\(A=\{a:(a,b)\in\Lambda\}\) is a set without multiplicities and satisfies
\begin{equation}
 \sum_{a\in A}
 \bigl|\langle f,e_a\rangle_{L^2(\mu_{1/2})}\bigr|^2
 =2\|f\|_{L^2(\mu_{1/2})}^2.
\end{equation}
It is a tight frame of bound \(c=2\) on an interval of length \(L=2\).
Corollary~\ref{cor:strip} gives strict pairwise separation greater than \(1/4\).
But Lemma~\ref{lem:packing}, with \(\delta=1/(cL)=1/4\), forces
\(A=a_0+\tfrac14\Z\), which contains pairs at distance exactly \(1/4\).
This is the required contradiction.
\end{proof}

\section{The parity dichotomy}\label{sec:parity}

We now combine the even-base construction with an \emph{external} obstruction to
prove the parity dichotomy announced in Corollary~\ref{cor:parity-intro}. The
even-base half is Theorem~\ref{thm:even}; the odd-base half relies on the recent
result that odd-base Cantor measures admit no exponential Fourier frame
\cite{DLT26,FSW26}. This dependence is explicit, and the odd-base nonexistence
theorem is not reproved here.

Recall from the introduction that, by Lai--Liu--Prince \cite{LLP21}, a Fourier
frame of an additive measure projects to Fourier frames of its components. The
next lemma is the symmetric, self-contained refinement we use: for a spectrum
the projections are additionally injective and the projected families are
exactly \(2\)-tight.

\begin{lemma}\label{lem:projection}
Let \(\mu\) be a continuous Borel probability measure on \(\R\) with
\(0\in\supp\mu\), and let
\[
 \frac12(\mu\times\delta_0+\delta_0\times\mu)
\]
be the associated additive measure. If it possesses a spectrum
\(\Lambda\subset\R^2\), then both coordinate projections of \(\Lambda\) are
injective and each projected exponential family is an exact \(2\)-tight Fourier
frame for \(L^2(\mu)\).
\end{lemma}

\begin{proof}
If distinct \((a,b),(a,d)\in\Lambda\) have the same first coordinate, then
orthogonality gives
\[
 \widehat\mu(d-b)=-1
\]
(up to conjugation, which is immaterial). But \(\widehat\mu(u)=-1\) is
impossible when \(0\in\supp\mu\): the average of the unimodular quantities
\(e^{-2\pi iux}\) can equal \(-1\) only if those quantities are \(-1\)
\(\mu\)-almost everywhere, and continuity then forces the same value at every
support point, including \(x=0\), where the value is \(1\). Thus the first
projection is injective; the same argument applies to the second projection.

For \(f\in L^2(\mu)\), apply Parseval in
\(L^2\bigl(\frac12(\mu\times\delta_0+\delta_0\times\mu)\bigr)\cong
L^2(\mu)\oplus L^2(\mu)\) to the vector \((f,0)\). Since the additive measure
puts weight \(\tfrac12\) on each axis,
\[
 \bigl\langle (f,0),e_{(a,b)}\bigr\rangle
 =\tfrac12\langle f,e_a\rangle_{L^2(\mu)},
\]
and Parseval gives
\[
 \|f\|^2
 =\frac12\sum_{(a,b)\in\Lambda}
   \bigl|\langle f,e_a\rangle_{L^2(\mu)}\bigr|^2.
\]
Injectivity turns the right side into an unweighted sum over the
first-coordinate projection. This is an exact \(2\)-tight Fourier frame.
The second-coordinate argument is identical.
\end{proof}

\begin{corollary}\label{cor:parity}
For the measures \eqref{eq:two-param} with integer \(R\ge2\),
\begin{equation}\label{eq:iff-even}
 \mu^{1/R}\text{ is spectral}
 \quad\Longleftrightarrow\quad
 R\text{ is even}.
\end{equation}
\end{corollary}

\begin{proof}
If \(R\) is even, apply Theorem~\ref{thm:even}.
If \(R\ge3\) is odd, the no-Fourier-frame theorem in \cite{DLT26} applies to
\(\mu_{1/R}\); rescaling the two equally weighted digits does not change
Fourier-frame existence. A spectrum for \(\mu^{1/R}\) would, by
Lemma~\ref{lem:projection}, give a Fourier frame for \(\mu_{1/R}\),
contradicting that theorem.
\end{proof}

\begin{corollary}\label{cor:odd-component}
Let \(R,S\ge2\) be integers. If at least one of \(R,S\) is odd, then
\(\mu_{1/R,1/S}\) is not spectral.
\end{corollary}

\begin{proof}
Suppose that \(\mu_{1/R,1/S}\) has a spectrum \(\Lambda\). The same argument as
in Lemma~\ref{lem:projection} shows that each coordinate projection of
\(\Lambda\) is injective: if two frequencies have the same first (respectively
second) coordinate, orthogonality would force
\[
 \widehat{\mu_{1/S}}(t)=-1
 \qquad
 \bigl(\text{respectively }\widehat{\mu_{1/R}}(t)=-1\bigr),
\]
impossible since \(0\) belongs to the support of both component measures. By
applying Parseval to functions supported on one axis, its projected frequency
set is a tight Fourier frame with bound \(2\) for that component. For any odd
base, this contradicts the no-Fourier-frame theorem of \cite{DLT26}.
\end{proof}

\begin{remark}
Dai's characterization of \emph{one-dimensional spectrality} alone would not
imply the odd-base half of \eqref{eq:iff-even}: failure to have an orthonormal
exponential basis is weaker than failure to have a Fourier frame. The latter is
what makes the projection argument decisive.
\end{remark}

\section{Concluding remarks}\label{sec:discussion}

\subsection{Example and open questions}
\begin{exam} Taking $k=2$ in Theorem~\ref{thm:even}, we have that the standard middle-fourth Cantor measure $\mu_{1/4}$ has spectrum \[ S_4=9\Gamma_4\cup(-3-9\Gamma_4), \]
 where $\Gamma_4=
      \left\{
        \sum_{j=0}^{m-1}\varepsilon_j4^j:
        m\ge 1,\ \varepsilon_j\in\{0,1\}\right\}$, while the planar measure $\mu^{1/4}$ has spectrum \[ \Lambda_4=\{(t,-t):t\in S_4\cup(S_4+3/4)\}. \] \end{exam}

We have worked throughout with the uncentered Bernoulli convolution
\eqref{eq:two-param}, generated by the digit set $\{0,2\}$. In the centered convention of \cite{Dai12}, the Bernoulli convolution is the
self-similar measure generated by the IFS
$\{x\mapsto r(x-1),\ x\mapsto r(x+1)\}$ with equal weights. Conjugating this
IFS by the translation $x\mapsto x+\frac{r}{1-r}$ yields the IFS
$\{x\mapsto rx,\ x\mapsto r(x+2)\}$ of \eqref{eq:ifs}; hence, up to a global
scaling, $\mu_r$ is a translate of the centered Bernoulli convolution. Translation is invisible to
one-dimensional spectrality---if $\tau_t\sigma$ denotes the translate of a
probability measure $\sigma$ on $\R^d$ by $t$, then
$\widehat{\tau_t\sigma}(\xi)=e^{-2\pi i\langle t,\xi\rangle}\widehat\sigma(\xi)$,
so the zero set of the Fourier transform and the Jorgensen--Pedersen spectral
function \eqref{eq:JP} are unchanged, and $\Lambda$ is a spectrum of $\sigma$
if and only if it is a spectrum of $\tau_t\sigma$. Hence our one-dimensional
statements match Dai's centered setting.

For the additive model \eqref{eq:two-param}, however, translation is not
invisible. A common translation of the whole planar measure moves both
coordinate axes together and is again harmless, but centering the two
components independently is not a single translation of the additive
measure; it is a genuine relative deformation of the two axes. This relative
translation is precisely the parameter governing the Lebesgue-type family in
the introduction, whose spectrality depends decisively on $t_1,t_2$
\cite{LLP21}. It is therefore natural to ask for which $t_1,t_2$ the
translated additive measures
\[
 \tfrac12\bigl((\tau_{t_1}\mu_r)\times\delta_0
       +\delta_0\times(\tau_{t_2}\mu_s)\bigr)
\]
are spectral. Our results concern only the uncentered case $t_1=t_2=0$; the
translation-dependence of spectrality for Bernoulli-type components remains
open.

Theorem~\ref{thm:unequal} provides an unconditional nonspectral
example with two different even bases. Nonstandard spectra of
one-dimensional Cantor measures have also been studied in
\cite{DHS09,DH15}. It does not imply that all
$\mu_{1/(2k_1),1/(2k_2)}$ with $k_1\ne k_2$ are nonspectral.
The arithmetic and unitary-gluing structures of two distinct singular
components require additional study.

The Bernoulli convolutions treated here are the two-digit instance of a large
class of self-similar spectral measures: those generated by a Hadamard triple. Let $q\ge2$ be an integer and let $D,C\subset\mathbb Z$ be finite sets with
$\#D=\#C$. If the matrix
\[
\frac{1}{\sqrt{\#D}}\bigl(e^{2\pi i dc/q}\bigr)_{d\in D,\,c\in C}
\]
is unitary, then $(q,D,C)$ is called a \emph{Hadamard triple}. For such a
triple, the self-similar measure $\mu_{q,D}$ generated by the maps
$x\mapsto q^{-1}(x+d)$, $d\in D$, with equal weights is spectral \cite{LW02}. Let $r\in \R$. It is easy to show that \begin{equation}\label{eq:scaldee} \text{$\Lambda$ is a spectrum of $\mu_{q,D}$ if and only if $r\Lambda$ is a spectrum of $\mu_{q,r^{-1}D}$.}
\end{equation}
Our component measure $\mu_{1/R}=\mu_{R,\{0,2\}}$ is obtained from $\mu_{R,\{0,1\}}$ by the dilation $x\mapsto 2x$, so it is a scaling of the above integer Hadamard triple model $(R,\{0,1\},\{0,R/2\})$.

It is natural to ask for which Hadamard triples $(q_1,D_1,C_1)$ and
$(q_2,D_2,C_2)$ the additive measure
\[
 \tfrac12\bigl(\mu_{q_1,D_1}\times\delta_0+\delta_0\times\mu_{q_2,D_2}\bigr)
\]
is spectral. Our results settle the symmetric two-digit case ($q_1=q_2=R$,
$D_1=D_2=\{0,2\}$) through the parity dichotomy, and exhibit the nonspectral
example $\mu_{1/2,1/4}$ for unequal contraction ratios $1/q_1\ne 1/q_2$. The general problem, including the cases $D_1\neq D_2$, $\#D_1=\#D_2>2$, and higher dimensions, remains open.

\subsection{Future directions}\label{subsec:program}

The results of this paper open a window on the spectral theory of the
additive model \eqref{eq:two-param}. For Cantor-type spectral measures, a
fairly complete structural picture has been built around three themes as mentioned in the introduction: the
spectral eigen-matrix problem (spectral eigenvalue problem in one dimension), the distribution of Beurling dimension of spectra, and the convergence of the associated Fourier
series. For the additive measures $\mu_{r,s}$, none of
these questions could even be formulated before, since the very spectrality
of the underlying planar measure was unknown. Theorem~\ref{thm:even} and
Corollary~\ref{cor:parity-intro} now supply exactly this missing starting
point, and thereby turn the following problems into well-posed ones.

(i) Spectral eigenvalue problem. Given a spectrum $\Lambda$ of
$\mu_{r,s}$, characterize the set of matrices $A\in M_2(\R)$ for which the dilated set
$A\Lambda$ is also a spectrum of $\mu_{r,s}$. This is completely open; the first case to
settle is the explicit spectrum $\Lambda_R$ of $\mu^{1/R}$ from Theorem~\ref{thm:even}.

(ii) Distribution of Beurling dimension of spectra. Recall that for a countable
set \(\Lambda\subset\R^d\), the {\it upper Beurling dimension} (or simply the {\it Beurling dimension}) of $\Lambda$ is defined by
\[
 \dim_{\mathrm{Be}}(\Lambda)
 =\limsup_{r\to\infty}\frac{\log N_\Lambda(r)}{\log r},
 \]
 where $N_\Lambda(r):=\sup_{x\in\R^d}\#\bigl(\Lambda\cap B(x,r)\bigr)$

 For the explicit spectrum $\Lambda_4$ of $\mu^{1/4}$, a direct count gives
$\dim_{\rm Be}(\Lambda_4)=\frac12$, which coincides with the entropy dimension of $\mu^{1/4}$. By Theorem 1.2 in \cite{Shi21}, the
Beurling dimension of any spectrum of a spectral measure is bounded above by
the entropy dimension of the measure; hence $\Lambda_4$ already attains the
maximal possible Beurling dimension. This focuses attention on the opposite
extreme: does $\mu^{1/4}$ admit a spectrum of Beurling dimension $0$? The
following proposition answers this question affirmatively.

\begin{proposition}\label{prop:zero-dim}
The measure \(\mu^{1/4}\) has a spectrum \(\Lambda_{\rm s}\) with
\(\dim_{\mathrm{Be}}(\Lambda_{\rm s})=0\).
\end{proposition}

We first recall a lemma.

\begin{lemma}[\cite{CKS08}]\label{incadd}
Let \(\Lambda_1,\Lambda_2\subseteq\R^d\) be countable sets. Then
\[
 \dim_{\mathrm{Be}}(\Lambda_1\cup\Lambda_2)
 =\max\bigl(\dim_{\mathrm{Be}}(\Lambda_1),\dim_{\mathrm{Be}}(\Lambda_2)\bigr).
\]
\end{lemma}

\begin{proof}[Proof of Proposition~\ref{prop:zero-dim}]
Consider the self-similar measure \(\mu_{4,\{0,6\}}\). By
\cite[Theorem~1.6]{AL23-0}, \(\mu_{4,\{0,6\}}\) has a spectrum
\(\Theta\subset\Z\) with \(\dim_{\mathrm{Be}}(\Theta)=0\). By
\eqref{eq:scaldee}, \(3\Theta\) is a spectrum of \(\mu_{1/4}\).

Set
\[
 T=(3\Theta)\cup\left(3\Theta+\frac34\right),
 \qquad
 \Lambda_{\rm s}=\{(t,-t):t\in T\}.
\]
We claim that \(\Lambda_{\rm s}\) is a spectrum of \(\mu^{1/4}\). The planar
computation in the proof of Theorem~\ref{thm:even} uses only the following
facts: a set \(S\) is a spectrum of \(\mu_{1/4}\), \(S\subset\frac32\Z\),
and \(2\cdot(3/4)\cdot1/3\equiv1/2\pmod1\). The set \(3\Theta\) satisfies these
properties: it is a spectrum of \(\mu_{1/4}\), it lies in
\(3\Z\subset\frac32\Z\), and \(2\cdot(3/4)=3/2\equiv1/2\pmod1\).
Therefore \eqref{eq:twotight}--\eqref{eq:Kcancel} apply verbatim and give
\[
 Q_{\Lambda_{\rm s}}^{\mu^{1/4}}(\xi,\eta)
 =\frac14\left(2+2+2\operatorname{Re}0\right)=1
 \qquad((\xi,\eta)\in\R^2).
\]
By Theorem~\ref{thm:JP}, \(\Lambda_{\rm s}\) is a spectrum of \(\mu^{1/4}\).

It remains to compute the Beurling dimension. By definition,
\(\dim_{\mathrm{Be}}\) is invariant under scaling and translation. Hence,
using \(\dim_{\mathrm{Be}}(\Theta)=0\),
\[
 \dim_{\mathrm{Be}}(3\Theta)=0,
 \qquad
 \dim_{\mathrm{Be}}\left(3\Theta+\frac34\right)=0.
\]
By Lemma~\ref{incadd},
\[
 \dim_{\mathrm{Be}}(T)
 =\max\left\{
 \dim_{\mathrm{Be}}(3\Theta),
 \dim_{\mathrm{Be}}\left(3\Theta+\frac34\right)
 \right\}=0.
\]

To bound \(\dim_{\mathrm{Be}}(\Lambda_{\rm s})\), let
\(\pi\colon\R^2\to\R\), \(\pi(x,y)=x\), be the first-coordinate projection.
Then \(\pi\) is \(1\)-Lipschitz, and its restriction to
\(\Lambda_{\rm s}=\{(t,-t):t\in T\}\) is a bijection onto \(T\). Moreover,
\[
 \pi\bigl(B_{\R^2}(z,r)\bigr)
 =\{x\in\R:\,|x-\pi(z)|<r\}
 =B_{\R}\bigl(\pi(z),r\bigr).
\]
Combining this with the bijectivity of \(\pi|_{\Lambda_{\rm s}}\), we obtain
\[
 \#\bigl(\Lambda_{\rm s}\cap B_{\R^2}(z,r)\bigr)
 \leq \#\bigl(T\cap B_{\R}\bigl(\pi(z),r\bigr)\bigr)
 \le N_T(r).
\]
Taking the supremum over \(z\in\R^2\) yields
\(N_{\Lambda_{\rm s}}(r)\le N_T(r)\), and therefore
\[
 \dim_{\mathrm{Be}}(\Lambda_{\rm s})
 \le\dim_{\mathrm{Be}}(T)=0.
\]
Since \(\dim_{\mathrm{Be}}(\Lambda_{\rm s})\ge0\), we obtain
\(\dim_{\mathrm{Be}}(\Lambda_{\rm s})=0\). This completes the proof.
\end{proof}

More generally, does the intermediate value property hold, that is, is every
value in $[0,\tfrac12]$ realized as the Beurling dimension of some spectrum
of $\mu^{1/4}$?

(iii) Convergence of Fourier series. For a spectrum $\Lambda$ of
$\mu_{r,s}$, every $f\in L^2(\mu_{r,s})$ has the formal expansion
$f\sim\sum_{\lambda\in\Lambda}\langle f,e_\lambda\rangle e_\lambda$; its
convergence is a new problem in this setting. Even the analogue of Carleson's theorem for the spectrum $\Lambda_R$ of $\mu^{1/R}$ is far from understood.

\section*{Acknowledgements}
The authors thank Ruxi Shi for raising the spectrality problem for additive measure with
Bernoulli convolutions in his talk at the ``Academic Symposium on Fractal
Spectral Measure Theory and Related Topics'' (Hunan Normal University,
May 15--18, 2026), which brought this question to our attention again. We also thank the organizers of the symposium. AI tools were used during the preparation of the paper to
assist with the elaboration of some proof details, checking calculations, and improving
the exposition. All mathematical statements and proofs were verified by the authors.
The authors take full responsibility for the paper.

\end{document}